\documentclass[12pt]{article}

\usepackage{amsmath,amsfonts, amssymb}

\begin{document}
\bibliographystyle{plain}
\pagenumbering{arabic}
\raggedbottom

\newtheorem{theorem}{Theorem}[section]
\newtheorem{lemma}[theorem]{Lemma}
\newtheorem{proposition}[theorem]{Proposition}
\newtheorem{corollary}[theorem]{Corollary}
\newtheorem{conjecture}[theorem]{Conjecture}
\newtheorem{definition}[theorem]{Definition}
\newtheorem{example}[theorem]{Example}
\newtheorem{condition}{Condition}
\newtheorem{main}{Theorem}
\setlength{\parskip}{\parsep}
\setlength{\parindent}{0pt}

\def \outlineby #1#2#3{\vbox{\hrule\hbox{\vrule\kern #1% 
\vbox{\kern #2 #3\kern #2}\kern #1\vrule}\hrule}}%
\def \endbox {\outlineby{4pt}{4pt}{}}%

\newenvironment{proof}
{\noindent{\bf Proof\ }}{{\hfill \endbox
}\par\vskip2\parsep}

\hfuzz12pt

\newcommand{\var}{\mbox{\rm Var}\;}
\newcommand{\supp}{\mbox{supp}}
\newcommand{\tends}{\rightarrow \infty}
\newcommand{\ep}{{\mathbb{E}}}
\newcommand{\pr}{{\mathbb{P}}}
\newcommand{\re}{{\mathbb{R}}}

\nocite{chen3} \nocite{utev}
\title{Convergence of the Poincar\'{e} constant}
\author{Oliver Johnson}
\date{\today}
\maketitle
\pagestyle{headings}

\makeatletter
\begin{abstract} 
The Poincar\'{e} constant $R_Y$ of a random variable $Y$ relates the 
$L^2(Y)$-norm
of a function $g$ and its derivative $g'$. Since $R_Y - \var(Y)$ 
is positive, with equality if and only if $Y$ is normal, it can be seen as a
distance from the normal distribution. In this paper
we establish the best possible 
rate of convergence of this distance in the Central
Limit Theorem. Furthermore, we show that $R_Y$ is finite for discrete
mixtures of 
normals, allowing us to add rates to the proof of the
Central Limit Theorem in the sense of
relative entropy.
\renewcommand{\thefootnote}{}
\footnote{{\bf Address:} Statistical Laboratory, 
Centre for Mathematical Sciences,
Wilberforce Road, Cambridge,  CB3 0WB, United Kingdom.  
Email: {\tt otj1000@cam.ac.uk}. }
\renewcommand{\thefootnote}{\arabic{footnote}}
\setcounter{footnote}{0}
\makeatother
\end{abstract}

\newlength{\dsl}
\newcommand{\doublesub}[2]{\settowidth{\dsl}{$\scriptstyle
#2$}\parbox{\dsl}{\scriptsize\centering { \normalsize $\scriptstyle #1$}
\\ {\normalsize $\scriptstyle #2$}}}

\section{Introduction and results}
Poincar\'{e} (or spectral gap) inequalities 
provide a relationship between $L^2$ norms on functions and their
derivatives.
\begin{definition}[Borovkov and Utev] \label{def:borov} 
Given a random variable $Y$, define
the Poincar\'{e} constant $R_Y$:
$$ R_Y = \sup_{g \in H_1(Y)} \frac{ \var g(Y)}{ \ep g'(Y)^2},$$
where $H_1(Y)$ is the space of absolutely continuous functions on the real 
line such that $\var g(Y) >0$ and $\ep g'(Y)^2 < \infty$.
\label{def:poincons} \end{definition}
$R_Y$ will not in general be finite, however it will be finite 
for the normal and other strongly unimodal
distributions (see for example Klaasen (1985)\nocite{klaasen}, Chernoff 
(1981)\nocite{chernoff}, Chen\nocite{chen2} (1982), 
Cacoullos\nocite{cacoullos} (1982), Nash\nocite{nash} (1958), 
Borovkov and Utev\nocite{borovkov2} (1984)). 

We will exploit various relationships between the Poincar\'{e} constant and
Fisher information:
\begin{definition} For a random variable $Y$ with smooth density $p$, define
the score function $\rho_Y(y) = p'(y)/p(y)$, and Fisher information $I(Y) =
\ep \rho(Y)^2$.
\end{definition}

Notice that for a given $Y$, if $g$ is a local maximum of
$\var g(Y)/(\ep g'(Y)^2)$ then for all functions $h$ and small $t$:
$$ \frac{\var(g +th)}{\ep(g'+th')^2}
\leq \frac{\var(g)}{\ep g^{'2}} = R_g,$$
so multiplying out, $0 \leq t^2 (R_g \ep h^{'2} - \ep h^2) +
2t (R_g \ep g' h' - \ep gh)$, which can only hold in an interval around
zero if:
\begin{equation} \label{eq:intpart} R_g \ep g' h' = \ep gh. \end{equation}
Integration by parts implies therefore that $g = - R_g(\rho_Y g' + g'') $,
so local maxima correspond to eigenfunctions of the Laplacian 
$D_Y g= (\rho_Y g' + g'')$, and the global maximum to the least strictly
negative eigenvalue (hence the alternative name of spectral gap inequality).

\begin{example}
The Poincar\'{e} constant can be infinite. For example, consider the
discrete random variable, where $\pr(X =1) = \pr(X= -1) =1/2$.
Then, we can choose $g$ such that $g'(-1) = g'(1)= \epsilon$, but
$g(-1) = -1$, $g(1)=1$, so that $\var g(X) =1$, but $\ep g'(X)^2 = 
\epsilon^2$. This 
argument will work for any discrete random variable, indeed any random
variable whose support is not an interval. \end{example}

However, our first main result shows that discrete random variables perturbed 
by small normals have a finite Poincar\'{e} constant:

\begin{theorem} \label{thm:discmix}
Consider $X$, a random variable with 
variance $\sigma^2$ taking a finite number of values 
with probabilities $p_1, p_2, \ldots p_n$ respectively,
and $Z_\tau$ an indepedent normal with mean zero variance $\tau$.
Then $Y_\tau = X + Z_{\tau}$  
satisfies a Poincar\'{e} inequality with constant $$R_{Y_{\tau}} = \tau 
\left(1 + \left( \frac{\sigma^2}{\tau \min_s p_s} \right)
\exp \left( \frac{\sigma^2}{\tau \min_s p_s} \right) \right).$$
\end{theorem}
\begin{proof} See Section \ref{sec:proofmix}. \end{proof}
Note that part 8 of Theorem 1.1 of Utev\nocite{utev} (1992) also shows that
$R_{Y_{\tau}}$ is finite. However, our bound has an explicit dependence on
$\sigma^2$, and so has independent interest.

In a paper by Johnson and Barron\nocite{johnson5} (2002), we 
show that finiteness of the Poincar\'{e} constant gives an explicit
rate of convergence
of relative entropy distance in the Central Limit Theorem. This is a strong
result, and implies convergence in $L^1$.

Theorems 2, 3 and 4 of Borovkov and Utev provide the following results:
\begin{lemma} \label{lem:borbas}
For the constant $R_X$ defined above:
\begin{enumerate}
\item{$R_{aX+b} = a^2 R_X$}
\item{If $X,Y$ are independent, then $R_{X+Y} \leq R_X + R_Y$}
\item{$R_{X} \geq \var(X)$, with equality if and only if $X$ is normal}
\item{If $R_X$ is finite then $\ep \exp (|X - \ep X|/12\sqrt{R_X}) \leq 2$,
so $X$ has moments of all orders.}
\item{If $R_{X_n}/\var(X_n) \rightarrow 1$, then $\ep w(X_n) \rightarrow
\ep w(Z)$ ,where $Z$ is normal, for any continuous $w$ with $|w(t)|
< \exp(c|t|)$, for sufficiently small $c$.} \end{enumerate} \end{lemma}
The first three properties are reminiscent of those of Fisher Information -- a 
subadditive relation holds and the minimising case characterises the 
normal distribution. In analogy with the approach to the Central Limit 
Theorem developed by Brown\nocite{brown} (1982) and Barron\nocite{barron}
(1986), we add an extra term into the subadditive relation
$R_{X+Y} \leq R_X + R_Y$, which is sandwiched as convergence occurs.
This gives us an answer to the question posed by Chen and Lou\nocite{chen3}
(1990),
of identifying the limit of the Poincar\'{e} constant in the Central Limit
Theorem. This was also answered by Utev\nocite{utev} (1992), though without the
explicit rate of convergence that we provide.
\begin{theorem} \label{thm:rconv}
Consider $X_1, X_2, \ldots $ IID, with $R=R_{X_i}$ and $I=I(X)$ finite. 
Defining $U_n = (X_1 + \ldots X_{n})/\sqrt{n \sigma^2}$,
then there exists a constant $C$, depending only on $I$ and $R$, such that
$R_{U_n} - 1 \leq C/n$. \end{theorem}
\begin{proof} See Section \ref{sec:proofconv}. \end{proof}
We can argue that this is the best possible rate, up to the choice of the
constant. Considering $g(x) = x^2-1$, we know that $R_X \geq (\ep X^4 -1)/4$.
Since $\ep U_{2n}^4 - 1 = (\ep U_n^4 - 1)/2 >0$, if $R_{U_n} - 1 = f(n)/n$
then $f(2^k)/2^k = R_{U_{2^k}} -1 \geq (\ep X^4 -1)/2^k$, so $f(2^k) \geq
\ep X^4 -1$ and hence does not tend to zero.

Since we can
 perturb random variables by adding small normals to ensure that the Fisher 
information is finite, we use this to
prove a strong form of the Central Limit Theorem.
\begin{theorem} If $X_1, X_2, \ldots $ are IID random variables with mean 0,
variance $\sigma^2$ and finite $R_X =R$, then $U_n = (X_1 + \ldots X_n)/
\sqrt{n \sigma^2}$ has the property that:
$$ \ep w(U_n) \rightarrow \ep w(Z),$$
where $Z$ is standard normal, for any continuous $w$ such that $w(t) \leq
\exp c|t|$, where $c < c_0 =  1/(12 \sqrt{R})$.\end{theorem}
\begin{proof} Given a random variable $U$, define $U^\tau \sim U + Z_\tau$.
Now, since $R_{U_n^\tau} \leq R_{U_n} + \tau \leq R + \tau$, by Lemma
\ref{lem:borbas}.4, taking $\tau_0$ such that $c = 1/(12 \sqrt{R + \tau_0})$:
$$ \ep \exp(c |U_n^\tau|) \leq 2,$$
for all $n$, if $\tau \leq \tau_0,c < c_0$.

Now, since $U_n^\tau$ have uniformly bounded Fisher information $I(U_n^\tau)
\leq I(Z_\tau) = 1/\tau$, and $R_{U_n^\tau} \leq R+\tau$, Theorem
\ref{thm:rconv} implies that $\ep w(U_n^\tau) \rightarrow \ep w(Z^\tau)$.

Hence, since:
\begin{eqnarray*}
\lefteqn{| \ep w(U_n) - \ep w(Z) |} \\
& \leq & | \ep w(U_n^\tau) - \ep w(U_n)| + 
| \ep w(U_n^\tau) - \ep w(Z^\tau)| + | \ep w(Z) - \ep w(Z^\tau)|,
\end{eqnarray*}
we need only show that given $\epsilon$, $|\ep w(U^\tau) - \ep w(U)| \leq
\epsilon$ for $\tau$ small enough.
This follows by uniform integrability arguments (see Theorem 25.12 of 
Billingsley\nocite{billingsley}), since $\ep |w(U^\tau)|^p \leq
\ep \exp cp|U_\tau| \leq 2$, for some small $p$, and since $w(U^\tau)$
converges weakly to $w(U)$.
\end{proof}

\section{\label{sec:proofmix} Finiteness of $R$ for mixtures of normals}

\begin{proof}{{\bf of Theorem \ref{thm:discmix}}}
Without loss of generality,
consider $X$ taking a finite number of values $a_1 > a_2 >
\ldots > a_n$ with probabilities $p_1, p_2, \ldots p_n$ respectively, 
where $\ep X = \sum_i p_i a_i =0$.

We introduce the `squared span' 
$M = \max (|a_1^2 - a^2_{2}|, |a_2^2 - a^2_{3}|,
\ldots |a_{n-1}^2 - a^2_{n}|, (a_1 - a_2)^2, (a_2 - a_3)^2, \ldots 
(a_{n-1} - a_n)^2 )$, and write $p$ for $\min_s p_s$.

By Theorem 1 of Borovkov and Utev, we need to check that for some $R$ and
all $x$, the density $f_\tau$ of $Y$ satisfies: 
\begin{equation} \label{eq:bornec}
 \int_x^{\infty} y f_\tau(y) dy \leq R f_\tau(x).\end{equation}

Since $f_\tau(y) = \sum_i p_i \phi_\tau(y-a_i)$, the LHS of 
Equation (\ref{eq:bornec}) becomes (defining $u_j = \sum_{i=1}^j p_i a_i$
and $a_{n+1} = - \infty$):
\begin{eqnarray*} 
\int_x^{\infty} y f_\tau(y) dy & = &
\sum_{i=1}^n p_i \int_x^{\infty} (y-a_i) \phi_\tau(y-a_i) dy 
+ \sum_{i=1}^n p_i a_i \int_x^{\infty} \phi_\tau(y-a_i) dy \\
& = & \sum_{i=1}^n p_i \tau \phi_\tau(x-a_i) dy
+  \sum_{i=1}^n p_i a_i \int_{x-a_i}^{\infty} \phi_\tau(y) dy \\
& = & \tau f_\tau(x) +
  \sum_{i=1}^n p_i a_i \sum_{j=i}^n \int_{x-a_j}^{x-a_{j+1}} \phi_\tau(y) dy \\
& = & \tau f_\tau(x) +
  \sum_{j=1}^{n-1} u_j \left( \int_{x-a_j}^{x-a_{j+1}} \phi_\tau(y) dy 
\right),
\end{eqnarray*}
since $u_n =0$, so for each interval $I_j = (x-a_j,x-a_{j+1})$
we need to consider bounds on $\min_{y \in I_j} y^2$.

 We write $r$ for the index such that $a_r \leq x < a_{r-1}$.

First, we consider $x \geq 0$, where
we can distinguish 3 cases: for $j <r$; $x-a_{j+1} < 0$, so for $y
\in I_j$:
$$ y^2 \geq (x - a_{j+1})^2 = (x- a_j)^2 + 2(a_j - a_{j+1})x + (a_{j+1}^2
- a_j^2) = (x - a_j)^2 - M.$$
For $j=r$; $y \in I_j$ means that:
$ y^2 \geq 0 \geq (x - a_j)^2 - (a_{j-1} - a_j)^2
\geq (x - a_j)^2 - M.$

For $j > r$; $x - a_j > 0$, so for $y \in I_j$:
$ y^2 \geq (x - a_{j})^2.$

Hence for all $j$, $\min_{y \in I_j} y^2 \geq (x-a_j)^2 - M$, so:
$$ \int_{x-a_j}^{x-a_{j+1}} \phi_{\tau}(y) dy
\leq (a_j - a_{j+1}) \max_{y \in I_j} \phi_{\tau}(y) \leq (a_j - a_{j+1}) 
\phi_{\tau}(x-a_j) 
\exp(M/2\tau).$$

In Lemma \ref{lem:techn}, we prove two technical results, that
$ u_j (a_j - a_{j+1}) \leq \sigma^2$, and that $M p \leq 2\sigma^2$.
This allows us to deduce that for $x \geq 0$:
\begin{eqnarray*}
\sum_{j=1}^{n-1} u_j \left( \int_{x-a_j}^{x-a_{j+1}} \phi_\tau(y) dy 
\right) & \leq &  \exp(M/2\tau) \sum_{j=1}^{n-1} u_j (a_j - a_{j+1}) 
\phi_{\tau}(x-a_j) \\
& \leq & \exp(M/2\tau) \sum_{j=1}^{n-1} \sigma^2  \phi_{\tau}(x-a_j)\\
& \leq & \exp(\sigma^2/\tau p) \sum_{j=1}^{n-1} \sigma^2 (p_j/p)
\phi_{\tau}(x-a_j),
\end{eqnarray*}
as required.

Similarly, for $x \leq 0$, we deduce that for $y \in I_j$, 
$ \min_{y \in I_j} y^2 \geq (x-a_{j+1})^2 -M$, and thus:
\begin{eqnarray*}
\sum_{j=1}^{n-1} u_j \left( \int_{x-a_j}^{x-a_{j+1}} \phi_\tau(y) dy 
\right) & \leq &  \exp(M/2\tau) \sum_{j=1}^{n-1} u_j (a_j - a_{j+1}) 
\phi_{\tau}(x-a_{j+1}) \\
& \leq & \exp(\sigma^2/\tau p) \sum_{j=1}^{n-1} \sigma^2 (p_{j+1}/p)
\phi_{\tau}(x-a_{j+1}).
\end{eqnarray*}\end{proof}

\begin{lemma} \label{lem:techn}
Using the notation above:
$ u_j (a_j - a_{j+1}) \leq \sigma^2$, and $M p \leq 2\sigma^2$.
\end{lemma}
\begin{proof} Note that $u_j = \sum_{i=1}^j p_i a_i = - \sum_{i=j+1}^n p_i 
a_i$.

For $a_{j+1} \geq 0$: 
$ u_j (a_j - a_{j+1}) \leq u_j a_j \leq \sum_{i=1}^j p_i a_i a_j
\leq \sum_{i=1}^j p_i a_i^2.$

For $a_{j} \leq 0$: 
$ u_j (a_j - a_{j+1}) \leq - u_j a_{j+1} \leq \sum_{i=j+1}^n p_i a_i a_{j+1}
\leq \sum_{i=j+1}^n p_i a_i^2.$

For $a_{j+1} \leq 0 \leq a_j$:
$ u_j (a_j - a_{j+1}) \leq \sum_{i=1}^j p_i a_i a_j u_j 
+ \sum_{i=j+1}^n p_i a_i a_{j+1}
\leq \sum_{i=1}^n p_i a_i^2.$

For the second part, we consider two cases, firstly where $M = a_s^2 - 
a_{s \pm 1}^2$. In this case:
$$ \sigma^2 = \sum_t p_t a_t^2 \geq p_s a_s^2 = p_s(a_s^2 - a_{s \pm 1}^2) + p_s
a_{s \pm 1}^2 \geq p M.$$
Alternatively, if $M = (a_s - a_{s+1})^2$ then:
$$ \sigma^2 = \sum_t p_t a_t^2 \geq p( a_s^2 + a_{s+1}^2)
\geq (p/2) (a_s - a_{s+1})^2 
\geq p M/2.$$ \end{proof}

\section{Convergence of the Poincar\'{e} constant} \label{sec:proofconv}
We establish an explicit rate of 
convergence of the Poincar\'{e} constant, using projection inequalities
similar to those in Johnson and Barron\nocite{johnson5} (2002)
\begin{lemma} \label{lem:rdec}
Given independent random variables $X,Y$ with Poincar\'{e}
constants $R_X,R_Y$, for any function $g$:
$$ \var g(X+Y) \leq (R_X + R_Y) \ep g'(X+Y)^2 - \frac{R_X}{R_X I(Y) + 1}
\var g'(X+Y),$$  and hence $R_{X+Y} \leq R_X + R_Y$. \end{lemma}
\begin{proof}
Without loss of generality, we can consider $g$ such that $\ep g(X+Y) =0$,
and define $h(u) = \ep_Y g(u+Y)$, which thus also has mean zero. Now:
\begin{eqnarray*}
\var g(X+Y) & = & \ep g^2(X+Y) \\
& = & \ep_X \left( \ep g^2(X+Y) | X \right)\\
& = & \ep_X \var (g(X+Y) | X) 
+ \ep_X \left( \ep g(X+Y) | X \right)^2 \\
& \leq & R_Y \ep_X ( \ep_Y g'(X+Y)^2 | X)
+ \ep h(X)^2 \\
& \leq & R_Y \ep g^{'2} (X+Y) 
+ R_X \ep h'(X)^2.
\end{eqnarray*}
To consider the second term, we 
use the score function $\rho_Y$ and define:
$$ f(x) = \ep_Y \left[ (g'(x+Y) - h'(x)) \rho_Y(Y) \right],$$
where by the Stein equation, $f(x) = - \ep_Y g''(x+Y) = - h''(x)$.
Further, by Cauchy-Schwarz:
$$ \ep h''(X)^2 = \ep f(X)^2 \leq I(Y) \ep (g'(X+Y) - h'(X))^2,$$
so that:
$$ \var h'(X) \leq R_X \ep  h''(X)^2 \leq R_X I(Y) \left( \ep g'(X+Y)^2 
- \ep h'(X)^2 \right),$$
and writing $\mu = \ep h'(X) = \ep g'(X+Y)$, we obtain:
$$ \ep h'(X)^2 (1+ R_X I(Y)) \leq R_X I(Y) \ep g'(X+Y)^2 + \mu^2,$$
which, rearranging, leads to:
$$ \ep h'(X)^2 \leq \ep g'(X+Y)^2 - \frac{\var g'(X+Y)}{R_X I(Y)+1}.$$
\end{proof}

Next we need a Lemma which again uses the idea that if $g'$ is nearly
constant, then $g$ is close to linear. We'd like to apply it to the optimal
$g$, which achieves the maximum in Definition \ref{def:borov}. However, 
rather than use compactness arguments to show such a function exists, we
can instead use a `good' $g$ instead.
\begin{lemma} \label{lem:gbd}
For any random variable $W$ with mean zero,  and any function $g$ such that
$R(t) = \var(g(W) +tW)/\ep(g'(W) +t )^2$ has a local maximum at $t=0$:
$$ \left( \frac{\var g(W)}{\ep g'(W)^2} - \var(W) \right) \leq
 3R_{W} \sqrt{\frac{\var g'(W)}{\ep g'(W)^2}}.$$
\end{lemma} \begin{proof}
Without loss assume that $\ep g(W) = 0$, and
write $\mu = \ep g'(W)$ and
$\delta_g = \var g'(W)/\ep g'(W)^2 = 1 - \mu^2/\ep g'(W)^2$ implies:
\begin{eqnarray*}
\lefteqn{\ep g^2(W)} \\
 &=& \ep (g(W) - \mu W)^2 + 2 \mu \ep W(g(W) - \mu W) + 
\mu^2 \ep W^2\\
& 
\leq & \ep (g(W) - \mu W)^2 + 2 |\mu| \sqrt{\var(W)} \sqrt{\ep (g(W) - \mu W)^2 }
+ \mu^2 \var(W) \\
& \leq & (\ep g'(W)^2) \left( R_W \left(\delta_g + 2 \sqrt{\delta_g
(1-\delta_g)} \right) + (1-\delta_g) \var(W) \right) \\
& \leq & (\ep g'(W)^2) \left( 3R_W \sqrt{\delta_g} + \var(W) \right) 
\end{eqnarray*} 
since $\ep (g(W) - \mu W)^2 \leq R_W \ep (g'(W) - \mu)^2 = R_W \var(g'(W))
\leq R_W \delta_g \ep g'(W)^2$, and since by Lemma \ref{lem:borbas}.3,
$\var(W) \leq R_W$. 
\end{proof}

Note, we can come up with tighter bounds: for example taking $h(W) =W$ in 
Equation (\ref{eq:intpart}),
$\ep W g(W) = R(0) \mu$. Hence
$\mu(R(0) - \var(W)) = \ep (g(W) - \mu W)W \leq \sqrt{\var(W)} 
\sqrt{ \ep (g(W) - \mu W)^2} \leq \sqrt{R_W \var(W) \var g'(W)}$.
This implies that:
$$ R(0) - \var(W) \leq R_W \sqrt{\frac{\delta_g}{1-\delta_g} }.$$
However, Lemma \ref{lem:gbd} is sufficient for our purposes.

\begin{proof}{{\bf of Theorem \ref{thm:rconv}}}
We consider convergence along the `powers of 2' subsequence $S_k = U_{2^k}$,
which implies convergence for the whole sequence by subadditivity.

For all $k$, $1 \leq R_{S_k} \leq R$, and $I_{S_k} \leq I$.
Taking $X =S_k/\sqrt{2}$ and $Y=S_{k}'/\sqrt{2}$ (an identical copy) in 
Lemma \ref{lem:rdec} implies (since $R_{S_k/\sqrt{2}} = R_{S_k}/2$ and
$I(S_k/\sqrt{2}) = 2I(S_k)$) that for any $g$:
$$\delta_g = \frac{\var g'(S_{k+1})}{\ep g'(S_{k+1})^2} \leq 
2(1/R_{S_k} + I(S_k)) \left( R_{S_k} -  
\frac{\var g(S_{k+1})}{\ep g'(S_{k+1})^2}\right).$$
Now, given $W = S_{k+1}$, we can find $h$ such that $\var h/\ep h^{'2} \geq 
\max(R - \epsilon, \var(W))$. Since $\var(h(W) +tW)/\ep(h'+t)^2$
tends to $\var(W)$ at $\pm \infty$, and has one maximum $t_0$ and one minimum,
we can find $g(W) = h(W)+ t_0 W$, which satisfies the conditions of 
Lemma \ref{lem:gbd}:
$$ \left( (R_{S_{k+1}} - \epsilon - 1)_+ \right)^2 \leq 
\left( \frac{ \var g(S_{k+1})}{\ep g'(S_{k+1})^2} - 1 \right)^2 
\leq 9 R^2_{S_{k+1}} \delta_g
\leq C(R_{S_k} - R_{S_{k+1}} + \epsilon),$$
where $x_+ = \max(x,0)$ and 
$C= 18R_{S_{k+1}}^2(1/R_{S_{k}}+I(S_{k})) \leq 18 R(IR+1)$.
That is, since $\epsilon$ is arbitrary,
\begin{equation} \label{eq:dynsys} (R_{S_{k+1}} - 1)^2 \leq C 
(R_{S_{k+1}} - R_{S_k}). \end{equation}
Note that 
since $R_{S_k}$ is decreasing and bounded below, successive differences
tend to zero, and thus $R_{S_k} \rightarrow 1$.

To obtain a rate, write $u_k = (R_{S_k} - 1)/C$, Equation (\ref{eq:dynsys}) 
gives $u_k (1 + u_k) \leq u_{k-1}$.
Since $u_k$ are decreasing:
$ \ldots u_{n+2}^2 \leq u_{n+1}^2 \leq u_n^2 \leq u_{n-1} - u_n$, and 
hence:
$ u_n^2 \leq u_{n-1} - u_n$, $ u_n^2 \leq u_{n-2} - u_{n-1}$,  
$ u_n^2 \leq u_{n-3} - u_{n-2}$ and so on. Summing, we obtain that for 
$m \leq n$:
$ (n-m) u_n^2 \leq u_m - u_n \leq u_m.$
Taking $n= 2^r$, $m= 2^{r-1}$ implies that:
$$ u_{2^r} \leq \sqrt{\frac{u_{2^{r-1}}}{2^{r-1}}}. $$
Repeating this $N$ times, we deduce that (since $u_k \leq u_1 \leq 1$):
$$ u_{2^r} \leq 2^{-\sum_{j=1}^N (r-j)/2^j}
= 2^{-r+2} 2^{(r-2-N)/2^N},$$
so if $N = r-2$, then $u_{2^r} \leq 4/2^r$, and we can `fill in
the gaps' by subadditivity, to show that $u_k \leq 16/k$ for all $k$.
\end{proof}

\section*{Acknowledgements}
The author is a Fellow of Christ's College, Cambridge, who helped support a 
trip
to Yale University during which many useful discussions with Andrew Barron
took place. Yurii Suhov of
Cambridge University and Alexander Holroyd of UCLA provided useful advice,
and the anonymous referee provided several extremely useful
improvements to the proofs.

\end{document}